\documentclass[12pt]{article}

\usepackage{diagrams}
\diagramstyle[l>=3em,scriptlabels,PostScript=dvips]

\usepackage{amsfonts}
\usepackage{amsmath}
\usepackage{amssymb}
\usepackage{amsthm}
\usepackage{graphicx}

\renewcommand{\H}{\mathbb{H}}
\newcommand{\F}{\mathbb{F}}

\newcommand{\Z}{\mathbb{Z}}

\newcommand{\Hom}{\operatorname{Hom}}
\newcommand{\Ext}{\operatorname{Ext}}
\newcommand{\Tor}{\operatorname{Tor}}
\newcommand{\HH}{HH} %% \operatorname{HH}}
\newcommand{\h}{\operatorname{H}}

\newtheorem{theorem}{Theorem}[section]
\newtheorem{lemma}[theorem]{Lemma}
\newtheorem{prop}[theorem]{Proposition}

\newtheorem{definition}{Definition}[section]
\newtheorem{remark}{Remark}

\begin{document}

\author{Dmitry Vaintrob}

\title{
The string topology BV algebra, Hochschild cohomology
and the Goldman bracket on surfaces
}

\date{}

\maketitle

\thispagestyle{empty}

\begin{abstract}

In 1999 Chas and Sullivan~\cite{cs} discovered that the homology
$\H_*(X)$ of the space of free loops on a closed oriented smooth manifold $X$
has a rich algebraic structure called \emph{string topology.}
They proved that $\H_*(X)$ is naturally a Batalin-Vilkovisky (BV) algebra.
There are several conjectures connecting the string topology BV algebra
with algebraic structures on the Hochschild cohomology
of algebras related to the manifold $X$, but none of them has been
verified for manifolds of dimension $n>1$.

In this work we study string topology in the case when $X$ is
aspherical (i.e. its homotopy groups $\pi_i(X)$ vanish for $i > 1$).
In this case the Hochschild cohomology Gerstenhaber algebra $HH^*(A)$
of the group algebra $A$ of the fundamental group of $X$ has a BV
structure. Our main result is a theorem establishing a natural
isomorphism between the Hochschild cohomology BV algebra $\HH^*(A)$ and
the string topology BV algebra $\H_*(X)$.
In particular, for a closed oriented surface $X$ of hyperbolic type
we obtain a complete description of the BV algebra operations
on $\H_*(X)$ and $HH^*(A)$ in terms of the Goldman
bracket~\cite{goldman} of loops on $X$.
The only manifolds for which the BV algebra structure on $\H_*(X)$
was known before were spheres~\cite{men} and complex Stiefel
manifolds~\cite{tam}.

Our proof is based on a combination of topological and algebraic
constructions allowing us to compute and compare multiplications
and BV operators on both $\H_*(X)$ and $\HH^*(A)$.
\end{abstract}

%% \newpage
%% \tableofcontents

\newpage

\section{Introduction}

For any topological space $X$ the set $G=\pi_1(X,x_0)$
of homotopy classes of based loops (continuous maps from the
circle $S^1$ to $X$ taking a fixed point $s_0$ in $S_1$ to
the base point $x_0$ in $X$) forms a group.
Forgetting the base points, we obtain the set $LX$ of
free loops, i.e.\ of all continuous maps $\phi:S^1 \to X$.
Homotopy classes of free loops are
indexed by conjugacy classes of the group $G$.
There is \emph{a priori} no multiplication or other algebraic
structure on this set. However, surprisingly,
if $X$ is an oriented two-dimensional manifold,
the vector space $L$ spanned by homotopy classes of free loops on $X$
has a natural Lie algebra structure. For two transversal loops
$\phi_1$ and $\phi_2$ their \emph{Goldman bracket} is defined by
\begin{equation}
\label{eq:goldm}
[\phi_1 ,\phi_2]=\sum_P s(P) \phi_1 *_P \phi_2,
\end{equation}
where the sum is taken over all intersections
$P$ of $\phi_1$ and $\phi_2$, the sign $s(P)\in\{\pm 1\}$ is
determined by the orientations ($s(P)=1$ if the area form
of the surface evaluated on the two tangent vectors is positive and
$s(P)=-1$ if it is negative) and the composition
$\phi_1 *_P \phi_2$ is the product of $\phi_1$ and $\phi_2$ as
elements of $\pi_1(X,P)$.

\begin{figure}[ht]
\label{fig1}
\begin{center}
\includegraphics[scale=0.2]{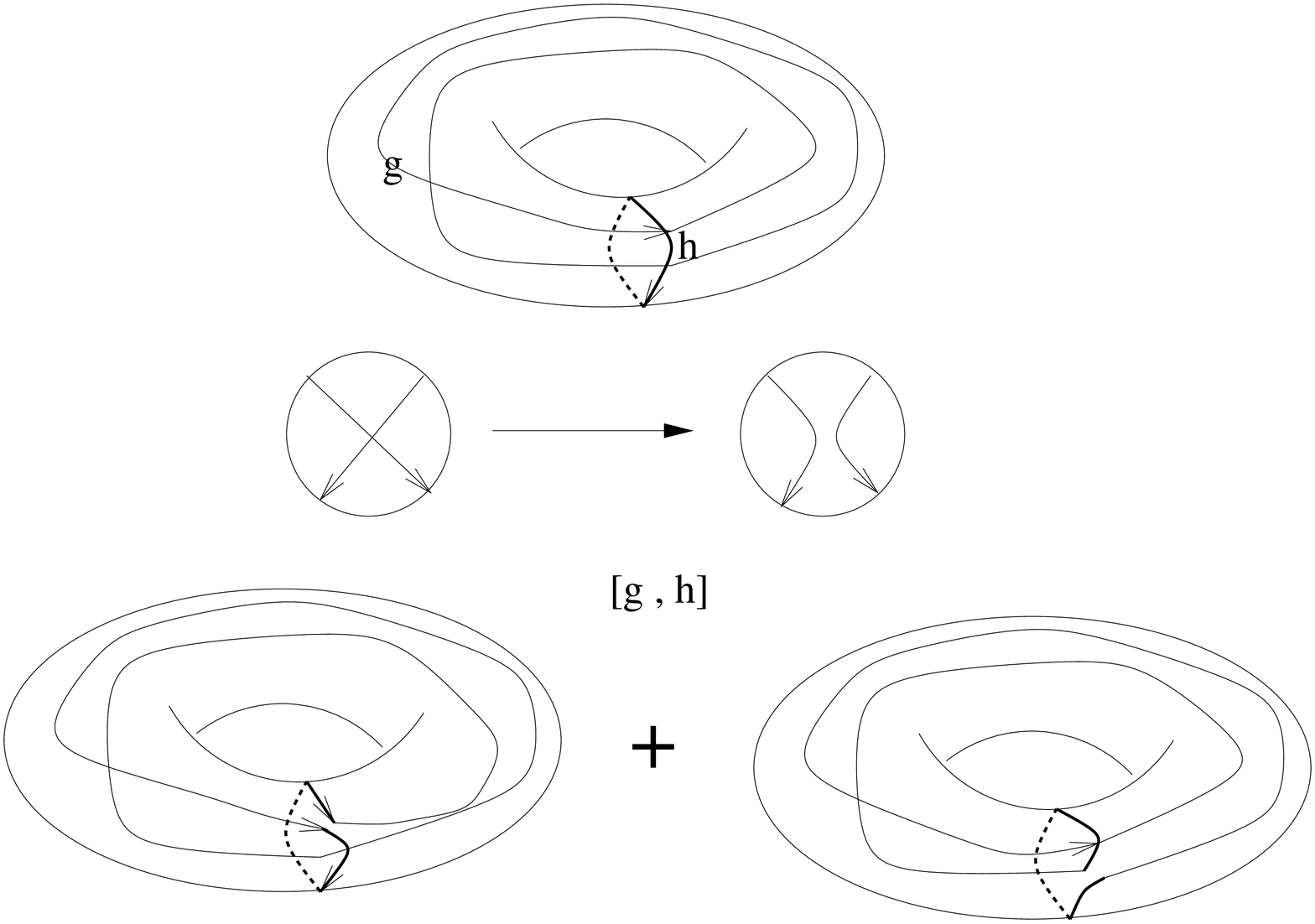}
\end{center}
\caption{The Goldman bracket of two loops on the torus}
\end{figure}

This Lie algebra structure was found by William Goldman
in~\cite{goldman} in his study of the symplectic structure on the moduli
space ${\Hom}(G,K)/K$ of representations of $G$ into a
compact Lie group $K$.

The Goldman Lie algebra was the starting point of the foundational work
by Chas and Sullivan~\cite{cs} on string topology.
Whereas Goldman considered only homotopy classes---or connected
components of the space of free loops $LX$ of a manifold $X$,
Chas and Sullivan studied algebraic structures on the
the total homology of $LX$. They called this space
$\H_*(X)=H_*(LX)$ the ``loop homology'' of $X$
and discovered that it has a natural multiplication, a Lie-type bracket,
and also a Batalin-Vilkovisky (BV) structure (see Section~\ref{sec:bv}).

String topology is a rapidly developing area of mathematical research
that connects classical algebraic topology with more recent developments
in mathematics influenced by theoretical physics and, in particular,
string theory and mirror symmetry (see e.g.~\cite{cv}).

The operations on the loop homology algebra of a manifold are very
difficult to compute. The complete BV algebra structure on $\H_*(X)$
has so far only been computed explicitly for
spheres~\cite{men} and complex Stiefel manifolds~\cite{tam} with
coefficients in an arbitrary ground field and for real projective
spaces~\cite{wes} with coefficients in $\F_2$. 

There are, however, several conjectures relating the loop
homology BV algebra with a more computable algebraic object,
the Hochschild cohomology of algebras related to $X$.

Hochschild cohomology $HH^*(A)$ of algebras was introduced by Hochschild
in 1945 as a tool for studying homological properties of algebras.
Recently it found important applications in other areas of mathematics
and also in theoretical physics.
In 1963 Gerstenhaber~\cite{ge} found that, in addition to the natural
cup product, $HH^*(A)$ has a Lie-type bracket. This bracket and the
cup product satisfy a compatibility relation and
make $HH^*(A)$ into a so-called Gerstenhaber algebra (see Section~\ref{sec:bv}).
Recently T. Tradler~\cite{tr1} showed that under certain assumptions
the Gerstenhaber structure on $HH^*(A)$ extends to a structure of a
Batalin-Vilkovisky algebra. In particular, this is the case
when $A$ is the algebra of singular cochains of a simply connected
closed oriented manifold (see~\cite{tr2}).
In~\cite{cj} Cohen and Jones established
that for a simply connected manifold $X$
there is an isomorphism
$$ F: {\H_*}(X) \to HH^*(C^*(X))$$
between the loop homology of $X$
and the Hochschild cohomology space of the
algebra $C^*(X)$ of singular cochains of $X$ and proved
that $F$ takes the loop product in ${\H_*}(X)$
to the cup product on Hochschild cohomology.
Since both ${\H_*}(X)$ and $HH^*(C^*(X))$ possess BV algebra
structure, it is conjectured that the Cohen-Jones identification $F$ is
an isomorphism of BV algebras. However, to the best of our knowledge,
the question whether $F$ respects the Gerstenhaber bracket or the BV operators 
remains open.

In this paper we algebraically compute the string topology BV algebra
for a large class of manifolds.
Namely, we show that for an aspherical smooth closed oriented manifold
$X$ of dimension $n$, its loop homology
BV algebra ${\H_*}(X)$ is isomorphic
to the Hochschild cohomology BV algebra $HH^*(A)$, where
$A$ is the group algebra $A=k[\pi_1(X)]$ of the fundamental
group of $X$ with the BV operator given by the dual
to the Connes operator $\kappa$ on the Hochschild homology $HH_*(A)$.

We first construct a vector space isomorphism between ${\H_*}(X)$
and $\HH_*(A)$ which sends the Chas-Sullivan BV operator $\Delta$
to the Connes operator $\kappa$. Then we construct a Poincar\'e duality
isomorphism $\tau: \HH_*(A) \to \HH^{n-*}(A)$ and show
that $\tau$ takes $\kappa$ to a BV operator for the Gerstenhaber
algebra on $\HH^{n-*}(A)$. We prove that the resulting
vector space isomorphism $\xi: {\H_*}(X) \to \HH^{n-*}(A)$
is an isomorphism of associative algebras and therefore
gives an isomorphism of BV algebras.

When $X$ is a closed oriented surface of genus $g\ge 2$
we obtain a complete description of the BV algebra operations
on $\H_*(X)$ and $HH^*(A)$ in terms of the Goldman
bracket~(\ref{eq:goldm}) of loops on $X$.

We hope that this result and our methods will provide an
insight for proving the Cohen-Jones conjecture that
$\H_*(X) \cong HH^*(C^*(X))$
for aspherical and other manifolds.

String topology on aspherical manifolds was also the subject
of the recent work~\cite{acg} by Abbaspour, Cohen and Gruher.
They described the loop homology product
in terms of a new
operation on the direct sum of group homologies of modules corresponding
to cosets of the fundamental group $G$ (not explicitly in terms of
Hochschild cohomology). However, they did not consider the Gerstenhaber
or the BV algebra structures.

\

This work originated as a project to compute
the Hochschild cohomology $HH^*(A)$ of the group
algebra $A$ of the fundamental group of a closed oriented hyperbolic
surface $X$ and express the Gerstenhaber structure
on $HH^*(A)$ in terms of the Goldman bracket~(\ref{eq:goldm}). 

This problem was motivated by the  following
result by Crawley-Boevey, Etingof and Ginzburg~\cite{cbeg}
about quiver algebras.

Let $P$ be the preprojective algebra of 
a hyperbolic (i.e.\  non-Dynkin and non-affine) quiver $Q$.
 The space $L=P/[P,P]$ has a natural Lie algebra structure given by the
 so-called \emph{necklace bracket}. Let
  $V_0$ be the vector space with basis given by the vertices of $Q$.
In~\cite{cbeg} it shown that 
$$HH^0(P)=V_0, \quad HH^1(P)=(L/V_0)\oplus \mathbf{k}, \quad HH^2(P)=L,$$
 and   Gerstenhaber algebra operations on $HH^*(A)$ can be expressed in terms
 of the necklace bracket.

{\bf Acknowledgements.}
I would like to thank the administrators and participants 
of the summer 2006 Research Science Institute at MIT where this work started.
I am grateful to my mentor Pavel Etingof for his guidance and
constant support. Finally, thanks are  to Mark Behrens, Allison Gilmore,
Christopher Michelich and Aaron Tievsky for their advice and
useful discussions. 

\section{Background and definitions}
\label{sec:background}

\subsection{Gerstenhaber and BV algebras}
\label{sec:bv}
Here we recall the definitions of Gerstenhaber and BV algebras
(see~\cite{get}).

\begin{definition}\em
A \emph{Gerstenhaber algebra} is a commutative graded algebra
$\displaystyle B^*=\bigoplus_{i\ge 0} B^i$ equipped with
a linear map
$$[\ ,\ ]: B^*\otimes B^*\to B^{*-1}$$
of degree $-1$ such that it
defines a super Lie algebra structure on the shifted space $B^{*-1}$ and,
for $a\in B^i$, the operator $[a,\cdot]$ is a degree $i-1$ derivation of
the product on $B^*$.

\end{definition}

\begin{definition}\em
\label{def:bv}
A \emph{Batalin-Vilkovisky (BV) algebra} is a commutative graded algebra $B^*$
with an operator $\Delta:B^*\to B^{*+1}$ (called the BV operator) such
that $\Delta\circ \Delta=0$ and the operation
$$[a,b]=\Delta(ab)-\Delta(a)b-(-1)^i a\Delta(b)$$
(where $a\in B^i$) defines a Gerstenhaber bracket on $B^*$.
\end{definition}

\subsection{Chas-Sullivan String Topology}
\label{sec:strtop}
Let $X$ be a closed oriented manifold of dimension $n$ and let $LX$ be
its space of free loops (maps $\gamma:S^1\to X$).

\begin{definition}
\em
The \emph{loop homology} ${\H_*}(X)$ of $X$ is 
the homology of its loop space, ${\H_*}(X)=H_*(LX)$.
Chas and Sullivan~\cite{cs} proved that
$\H_*(X)$ forms a BV algebra (see below).
We will call it the \emph{string topology BV algebra}.
\end{definition}

The BV algebra operations on $\H_*(X)$ are as defined as follows.

The \emph{loop product}
\begin{equation}
\label{eq:loop}
\H_i(X)\otimes \H_j(X)\to \H_{i+j-n}(X)
\end{equation}
is given by the combination of the intersection product on $X$,
\begin{equation}
\label{eq:intprod}
\cap:H_i(X)\otimes H_j(X)\to H_{i+j-n}(X),
\end{equation}
with the composition of loops
$LX\times_X LX \to LX.$ % make displayed later
Here the projection $LX\to X$ is given by the map $ev:\gamma\mapsto
\gamma(s_0)$ where $s_0\in S^1$ is the base point of the circle.

The \emph{BV operator}
\begin{equation}
\label{eq:delta}
\Delta:{\H_*}\to \H_{*+1}
\end{equation}
is determined by the natural $S^1$ action
$S^1\times LS \to LS$ % make displayed later
induced by the action of $S^1$ on itself. Note that the operator
$\Delta$ on $\H_*(LX)$ is defined even if $X$ is not a manifold.

\subsection{Hochschild homology and cohomology}

We recall here the definitions of Hochschild homology and
cohomology (see e.g.~\cite{w}.)

\begin{definition}
\em
Let $M$ be a bimodule over an algebra $A$.
The \emph{Hochschild homology} $HH_*$ and \emph{cohomology} $HH^*$
of $A$ with coefficients in $M$
are given by
\begin{equation}
\label{eq:homol}
HH_i(A,M)=\Tor_i(A,M) \quad \mathrm{and} \quad
HH^i(A,M)=\Ext^i(A,M),
\end{equation}
where the functors $\Tor_*$ and $\Ext^*$ are
taken in the category of $A$-bimodules.
\end{definition}

When $M=A$, the graded spaces $HH_*(A)=HH_*(A,A)$ and
$HH^*(A)=HH^*(A,A)$ have several natural algebraic structures.
Hochschild homology has a natural operator
\begin{equation}
\label{eq:connes}
\kappa:HH_i(A)\to HH_{i+1}(A)
\end{equation}
called the Connes operator. It comes from the cyclic structure on the
Hochschild complex of $A$ (see \cite{loday}).

Hochschild cohomology has two product structures, a cup product
\begin{equation}
\label{eq:cup}
\cup:HH^i(A)\otimes HH^j(A)\to HH^{i+j}(A)
\end{equation}
(which is associative and
graded commutative), and a Gerstenhaber Lie-type bracket~\cite{ge,loday}
\begin{equation}
\label{eq:gerst}
[\ ,\ ]:HH^i(A)\otimes HH^j(A)\to HH^{i+j-1}(A).
\end{equation}
Together these two structures make $HH^*$ into a Gerstenhaber algebra.

Van den Bergh \cite{van}
showed that for certain algebras $A$, there
exists a non-negative $n$ and a Poincar\'e duality isomorphism between
Hochschild homology and cohomology
\begin{equation}
\tau:HH_i\to HH^{n-i}.
\label{eq:dual}
\end{equation}
In a recent preprint \cite{gi2}, Ginzburg proved further that 
this isomorphism takes the Connes operator $\kappa$ on $HH_*$ to a BV 
operator $\lambda=\tau\circ\kappa\circ\tau^{-1}$ on $HH^*$ compatible with the
Gerstenhaber algebra structure.

\section{Statement of results}
\label{sec:res}

\subsection{String topology BV algebra of aspherical manifolds}

The following theorem is the main result of this paper.

\begin{theorem}
\label{thm:main}
Let $X$ be a closed connected 
oriented aspherical manifold of dimension $n$
and let $A$ be the group algebra of the fundamental group $G=\pi_1(X)$.
There exists an isomorphism of BV algebras
$$
\xi: \H_*(X) \to HH^{n-*}(A)
$$
where $\H_*(X)=H_*(LX)$ is the Chas-Sullivan loop homology
BV algebra and $HH^*(A)$ is the Hochschild cohomology algebra
equipped with the BV operator $\lambda$.
\end{theorem}

The proof of this theorem is given in the next section.

\subsection{Hochshild cohomology and Goldman bracket on surfaces}

Theorem~\ref{thm:main}  allows us to compute the Hochschild
cohomology BV algebra $HH^*(A)$ of the group algebra of 
the fundamental group  and  the string topology BV algebra $\H_*(X)$
for a hyperbolic surface $X$ in terms of the  Goldman Lie
bracket~(\ref{eq:goldm}) on $X$.

\begin{theorem} \label{thg}
Let $X$ be a compact oriented surface of genus $g>1$. Let $A=\mathbf{k}[G]$
be the group algebra of its fundamental group $G$ and let $L=H_0(LX)$
be the space generated by homotopy classes of free loops equipped with the
Goldman Lie bracket.

\noindent
(i) The Hochschild cohomology graded space of $A$ is naturally isomorphic to
the loop homology space, namely
\begin{equation}\label{eq:gold_hoch}
HH^i(A)\cong H_{2-i}(LX) \ \mathrm{for} \
0\le i \le 2, \quad \mathrm{and} \
HH^i(A)=0 \ \mathrm{for} \ i>2.
\end{equation}
\noindent
(ii) 
Under the above identification, the cup product on the algebra
$HH^*(A)$ coincides with the Chas-Sullivan loop homology product
(of degree $-2$) on $\H_*(X)$.

\noindent
(iii) Under the identification~(\ref{eq:gold_hoch}) the Gerstenhaber bracket
on $HH^*(A)$ becomes the Gerstenhaber bracket on the string topology
$\H_{2-*}(X)$.

\noindent
(iv) The Gerstenhaber algebra $HH^*(A)$ has a Batalin-Vilkovisky
structure given by the operator $\lambda$ of degree $-1$ that
corresponds to the Chas-Sullivan string topology BV operator $\Delta$
on $\H_*(X)$ (of degree $+1$).

\noindent
(v) The non-trivial Hochschild cohomology groups of $A$ are given
by
$$
HH^0(A)=\mathbf{k}, \ HH^1(A)=H_1(X,\mathbf{k})\oplus L/\mathbf{k}\gamma_0, \ HH^2(A)=L,
$$
where
$\gamma_0\in L$ is the class of the trivial loop.

\noindent
(vi) Non-trivial cup products on $HH^*(A)$ exist only on $HH^1(A)$
and can be expressed in terms of the Goldman bracket on $X$ as follows
$$ (\alpha,\gamma)\cdot (\alpha',\gamma')=
\langle \alpha,\alpha' \rangle\gamma_0+
\langle\alpha',\gamma\rangle\gamma +
\langle\alpha,\gamma'\rangle\gamma' +
[\gamma,\gamma'],$$
where $\langle\cdot,\cdot\rangle$ is the intersection pairing on
$H_1(X,\mathbf{k})$ and $[\cdot,\cdot]$ is the Goldman bracket on $L$.

\noindent
(vii) The BV operator $\lambda$ is equal to $0$ on $HH^1(A)$ and
is induced by the projection $L\to L/\mathbf{k}\gamma_0$ on $HH^2(A)=L$.
\end{theorem}
\proof

Parts (i)-(iv) follow from Theorem~\ref{thm:main}.

\noindent To show (v) and (vii), we see that $HH^*(A)$ as a space decomposes into
the sum over conjugacy classes $C$ of group cohomology
$$HH^*(A)=\bigoplus_C H^*(G,\mathbf{k} C)$$
where $\mathbf{k} C$ is the vector space spanned by the elements of $C$ on which 
$G$ acts by conjugation. We can
compute $H^*(G,\mathbf{k} C)$ using the fact that for any $g\in G, g\ne 1$, its
centralizer $Z(g)\cong \Z$. 
Since conjugacy classes of $G$ correspond to homotopy classes of free loops,
this gives the desired result.

\noindent Finally, to prove (vi) we use the formula
$[\gamma,\gamma']=\Delta(\gamma)\cdot\Delta(\gamma')$ 
which expresses the Goldman bracket of $\gamma, \gamma'\in L\cong\H_0(X)$
in terms of string topology operations.
 \qed

\section{Proof of Theorem~\ref{thm:main}}
\label{proofs}

\subsection{Notation and conventions}

We will be using the following notation and conventions.

\begin{itemize}
\item
$\mathbf{k}$ is a field of characteristic zero.

\item $X$ is a closed oriented smooth aspherical manifold of dimension
$n$, a base point $x_0$ and the fundamental group $G=\pi_1(X,x_0)$.

\item $A=\mathbf{k}[G]$ is the group algebra of $G$.

\item The notation $_A Y$ means that $Y$ is a left $A$-module.
By $Y_A$ we denote a right $A$-module and
by $_A Y_A$ an $A-A$ bimodule structure on $Y$.
Note that since $A=k[G]$ we have $A\cong A^{op}$ and therefore
each left $A$-module is canonically a right $A$-module.
In particular, this implies that $A\otimes A$-modules
can be viewed as $A-A$ bimodules.

\item
By $Y_*$ or $Y^*$ we denote a graded vector space or a chain complex,
and $Y_{n-*}, Y_{*-1}$ etc.\ denote the same thing with shifted grading.

\item By a map $f: X_*\to Y_*$ or $f:X_*\to Y_{n-*}$ etc.,
we always mean a homomorphism of complexes or graded spaces.

\item
Tensor product $\otimes$ is taken over $\mathbf{k}$ by default, $\otimes_A$ means
the tensor product over $A$,
and 
$\otimes_{A \otimes A}$ denotes the tensor product in the category of
$A-A$ bimodules.

\item 
By $C_*(X)$ we denote the singular chains of a
topological space $X$. For a cell decomposition  
 $T$ of $X$, we denote by $C_*(T)$ the corresponding chain
complex. 

\end{itemize}

\subsection{Construction of resolutions}

In this section we construct several resolutions of the $A$-module
$\mathbf{k}$ and of the  bimodule $_A A_A$.

Let $\tilde{X}$ be the universal covering space of $X$. Since $X$ is
aspherical, $\tilde{X}$ is contractible.
We define $\tilde{R}_*=C_*(\tilde{X}_i)$,
the $i$-dimensional singular chains on $\tilde{X}$.

\begin{prop}\label{prop:res1}
$\tilde{R}_*$ is a projective $A$-module resolution of $_A\mathbf{k}$.
\end{prop}

\proof
Since $G$ acts freely on $\tilde{X}$, $G$ also acts freely on
$C_*(\tilde{X})$ and so $\tilde{R}_*$ is a complex of projective
(indeed of free) $A$-modules. Since $H_*(\tilde{X})$ is equal to $\mathbf{k}$
concentrated in degree zero, $\tilde{R}$ is a projective resolution of
$_A\mathbf{k}$.
\qed

We will also use $A$-module resolutions of $\mathbf{k}$ using cellular chains.
Let $T_*=T_0\cup T_1\cup\ldots\cup T_n$ be a cellular decomposition of
$X$ (where $T_1$ are the 1-simplices, etc.) By the homotopy lifting
property of covering spaces, $T$ can be lifted to a cell decomposition
$\tilde{T}$ of $\tilde{X}$.

\begin{prop}
$C_*(\tilde{T})$ is a projective $A$-module resolution of $\mathbf{k}$.
\end{prop}
\proof
This is proven analogously to \ref{prop:res1}. \qed

We will work with two particular cell decompositions.

\begin{definition}
\em
Let $T$ be a triangulation of $X$ (one exists because $X$ is a smooth
manifold), and let $T'$ be the dual cellular decomposition.

We lift these decompositions to decompositions $\tilde{T}$ and
$\tilde{T}'$ respectively of $\tilde{X}$. We denote the resolutions
$R_*=C_*(\tilde{X})$ and $R'_*=C_*(\tilde{X}')$.
\end{definition}

We perform the same constructions on the fiber product
$Q=\tilde{X}\times_X \tilde{X}$
of the universal covering of $X$ with itself.
\begin{definition}
\em
Let $T^b$ be the simplicial complex given by the barycentric
subdivision of $T$. We lift this decomposition using the covering
homotopy property to get a cellular decomposition $T^b_Q$ of $Q$. We
define
$$W_*=C_*(T_Q^b)
\mathrm{\ \ and\ \ } \tilde{W}_*=C_*(Q).$$
\end{definition}

The $G$-action on $\tilde{S}$ defines a $G\times G$ action on $Q$.
This makes $Q$ an $A$-bimodule.
\begin{lemma}
The complexes $W_*$ and $\tilde{W}_*$ are projective $A$-bimodule
resolutions of  $_A A_A$.
\end{lemma}
\proof The modules $W_*$ and $\tilde{W}_*$ have free $G\times G$ action
and are therefore projective. All connected components of $Q$ are
contractible, so $H_*(Q)$ is concentrated in degree $0$. It is a
well-known result that the connected components of $Q$ are indexed by
elements of $G$ and that $H_0(Q)=A$ with canonical $A$-biaction. This
proves the proposition. \qed

\begin{remark}
We will later consider every point $\tilde{x}\in \tilde{X}$ as a point
$x\in X$ and homotopy class of paths $\epsilon$ from the base point
$x_0$ to $x$. Analogously, every point $q\in Q$ can be considered as
a pair of homotopy classes $\epsilon_1$ and $\epsilon_2$ from $x_0$ to
$x$. The connected component of $Q$ containing $q$ is then indexed by
the path product $\epsilon_2^{-1}\epsilon_1$.
\end{remark}

\subsection{Construction of the isomorphism $\rho$}
We will define a vector space isomorphism $\rho: \H_*\to HH_*$
which takes the BV operator $\Delta$ on string homology to the Connes
operator $\kappa$ on $HH_*$.

Let $_A{N}$ be the $A=k[G]$-module whose underlying space is $A$ and
with $G$-action defined by conjugation, $g.a=gag^{-1}$ for $g\in G$ and
$a\in A$.

We will use the following standard fact about Hochschild (co)homology of
group algebras.

\begin{prop}
For a group algebra $A=\mathbf{k}[G]$, its Hochschild homology and cohomology is
isomorphic to the group homology and cohomology of $G$ with coefficients
in $N$:

$$ HH_*(A)=\h_*(G,{N}) \ \ \mathrm{and} \ \ HH^*(A)=\h^*(G,{N}). $$
\end{prop}
A proof can be found in e.g.~\cite{loday} or \cite{w}.
\qed

Let $G$ be any (discrete) group and let $X=K(G,1)=BG$ (i.e.~$X$ is a
connected topological space with $\pi_1(X)=G$, $\pi_i(X)=0$ for $i\ge 2$).

Even if $X$ is not a manifold, the operator $\Delta:\h_*(LX)\to
\h_{*+1}(LX)$ is well-defined.

We will use the following known result.
\begin{theorem}[\cite{loday} Corollary 7.3.13]
There is an isomorphism of vector spaces,
\begin{equation}
\label{eq:bg}
\rho:\h_*(LX)\to HH_*(k[G])
\end{equation}
which takes the operator $\Delta$ to the Connes operator
$\kappa$ on $HH_*$.
\end{theorem}

Loday constructs this isomorphism in terms of the geometric realization
of the cyclic bar construction and its covering of the geometric
realization of the regular bar construction.

In terms of the resolution $\tilde{R}_*$, the map $\rho$ can be computed
as follows.

Let $\sigma \in C_i(LX)$ be an $i$-simplex, $\sigma:\Delta^i\to LX$.
Composition with the map $\mathrm{ev}_{s_0}:LX\to X$ gives us a simplex
$\sigma_0:\Delta^i\to X$. By the homotopy covering theorem, we can
(not canonically) choose a lifting of the map $\sigma_0$ to a map
$\tilde{\sigma}_0: \Delta^i\to \tilde{X}.$
Let $p\in \Delta^i$ be a point in the simplex. The point gives a loop
$\gamma=\sigma(p)$ and a point $\tilde{m}=\tilde{\sigma}(p)\in
\tilde{X}$. The point $\tilde{m}$ corresponds to a homotopy class
$\epsilon$ of paths from the base point $m_0$ to the point
$m=\sigma_0(p)$. The loop $\gamma$ represents a homotopy class $g_m$
of loops with base point $m$. Then the conjugate
$\epsilon^{-1}g_m\epsilon$ is a homotopy class of loops based at
$m_0$ and gives an element of the fundamental group $g\in G$. Note
that the element $g$ is independent of choice of $p\in \Delta$.

\begin{lemma}
The element
$g\otimes\tilde{\sigma}\in N\otimes_G \tilde{R}$
is independent of the choice of the lifting
$\tilde{\sigma}$ and therefore well-defined.
\end{lemma}
\proof Let
$\tilde{\sigma}'$ be a different lifting, and let $g'$ be the element of
$G$ which we get by the above construction. There is an element $h\in
G$ such that $\tilde{\sigma}'=h \tilde{\sigma}$ and this makes
$g=h^{-1}gh$ (the action of $h^{-1}$ on $g\in N$). The actions of $h$ and
$h^{-1}$ get canceled after taking tensor product over $A$, so
$g\otimes_G\tilde{\sigma}=g'\otimes_G\tilde{\sigma}'.$ \qed

We define $\rho_0:C_*(LX)\to N\otimes_G \tilde{R}$ by
$\rho_0(\sigma)=g\otimes\tilde{\sigma}.$
\begin{lemma}
  The map $\rho:\H_*\to HH_*$ induced by $\rho_0$ on homology coincides
  with the isomorphism given in~\cite[Corollary 7.3.13]{loday}.

\end{lemma}
\proof The geometric realization $|B.G|$ of the bar construction of $G$
is a $K(G,1)$ space and therefore homotopic to $X$. Any homotopy
equivalence $\psi: |B.G|\to X$ gives a quasiisomorphism of complexes
which identifies our construction with Loday's on the level of homology.
\qed

\subsection{Intersection product}

Now we will construct a non-commutative analogue of the intersection
product
\begin{equation}
  \label{eq:mu}
\mu:R_i\otimes{R}_j'\to W_{i+j-n}.
\end{equation}

Let $\sigma,\sigma'$ be simplices in
${T}_i,{T}'_j$ respectively and let $\sigma_0$ and
$\sigma_0'$ be their images in $T$ and $T'$. If $\sigma_0$ and
$\sigma_0'$ do not intersect, we set $\mu(\sigma\otimes \sigma')=0$.
Otherwise, from the definition of the dual complex, we know that
$\sigma_0$ and $\sigma_0'$ intersect in exactly one cell of the
barycentric subdivision, $\sigma_0^b\in T^b_{i_j-n}$. For every point
$p$ of $\sigma_0^b$, the cell $\sigma$ gives one point in its
$\tilde{X}$ fiber and $\sigma'$ gives another. This naturally gives
us a lifting of $\sigma_0^b$ to a cell
$\sigma^b$ of $\tilde{X}\times_X\tilde{X}$. We have
$\sigma^b\in W^b$ and we define
$\mu(\sigma\otimes\sigma')=\sigma^b.$
We extend $\mu$ to all of ${R}\otimes{R}'$ by linearity.

\begin{lemma}
The map $\mu$ respects the $A$-biaction, i.e.\ for any $a, b\in A$
$$\mu(a\sigma_i\otimes b\sigma_j)=(a\otimes
b)(\mu(\sigma\otimes\sigma')).$$
\end{lemma}
\proof This can be verified directly from the definition of $\mu$.
\qed

\

Using the map $\mu$, we will define a product
 $$\beta_0:C\otimes C'\to A _A\otimes_A W_{i+j-n}$$
such that the induced map $\beta:HH_i\otimes HH_j\to HH_{i+j-n}$ 
makes the following diagram commute:

\begin{equation}
\label{two-tier-diagram}
\begin{diagram}
\H_i \otimes \H_j & \rTo^\bullet & \H_{i+j-n}\\
\dTo^{\rho\otimes\rho} & & \dTo_\rho \\
HH_i\otimes HH_j & \rTo^\beta & HH_{i+j-n}\ .\\
\end{diagram}
\end{equation}

\begin{lemma}
\label{constr:beta}
Let $g\otimes_G\sigma\in C_i$ and $g'\otimes_G\sigma'\in C'_j$.
The equation
\begin{equation}
  \label{eq:prop}
\beta_0(g\otimes_G\sigma\otimes g'\otimes_G\sigma')=
g(\mu[\sigma\otimes\sigma']^{-1}g'\otimes_A^A(\mu(\sigma\otimes\sigma'))
\end{equation}
gives a well-defined map
 $\beta_0:C\otimes C'\to A _A\otimes_A W_{i+j-n}$
where $[\mu(\sigma\otimes\sigma')]\in G$ is the element of $G$
corresponding to the connected component of
$Q=\tilde{S}\times_S\tilde{S}$ which contains the cell
$\mu(\sigma\otimes\sigma')$.
\end{lemma}

\proof
This can be verified by a direct computation.
\qed

\begin{lemma}
The map $\rho$ takes the product $\cdot$ on string topology to the
operation $\beta$ on $HH_*$.
\end{lemma}
\proof
Let $\gamma\in \H_i, \gamma'\in \H_j$. Let $c$ be a representative of
$\rho(\gamma)$ in the complex $C_i$ and let $c'$ be a representative of
$\rho(\gamma')$ in $C'_j$.

It is possible to choose a representative $\gamma_0\in C_i(LX)$ of
$\gamma$ such that $\rho_0(\gamma_0)=c_0$
and similarly a representative $\gamma_0'\in C_j(LX)$
of $\gamma'$ such that
$\rho_0(\gamma_0')=c_0'$.

Since all intersections between cells of $T$ and cells of $T'$
are transversal, we can explicitly construct a representative of $\gamma
\cdot \gamma'$ by multiplying together all pairs of loops which map $s_0\in
S^1$ to the same point in $X$. We apply this operation to $\gamma_0$ and
$\gamma'_0$ to get a new chain of loops, $\gamma''_0$. From our
construction of $\beta_0$, we see that
$\rho_{X_0}(\gamma''_0)=\beta_0(\gamma\otimes \gamma').$
It follows that the two products, $\cdot$ and $\beta$, coincide on homology.
\qed

\subsection{Hochschild Poincar\'e Duality}

The pairing 
$$\pi_0=\alpha \circ \mu: {R}_i\otimes {R}'_{n-i}\to A$$
 has two adjoint maps
$$\iota:R_i\to \Hom_A(R'_{n-i}, A) \ \ \mathrm{and} \ \ \iota':R'_{n-i}\to
\Hom_A(R_{i}, A).$$

\begin{lemma}
The maps $\iota:{R}_i\to \Hom_A({R}'_{n-i}, A)$ and
$\iota':{R}'_{n-i}\to \Hom_A({R}_{i}, A)$ are isomorphism of
$A$-modules.
\end{lemma}
\proof $R$ is a free $A$-module. Therefore $\Hom(R,A)$ has a basis (as a
vector space) of maps $f_\sigma$ for $\sigma\in T_i$ with
$f_\sigma(\sigma)=1$ and $f_\sigma(\sigma_1)=0$ for a cell $\sigma_1\ne
g \sigma$ for some $g\in G$.
We see that for a cell $\sigma'\in T'_{n-i}$ the map $\iota':{R}'_i\to
\Hom_A({R}_{n-i}, A)$ takes $\sigma'$ to $f_\sigma$ where $\sigma$ is the
unique cell which intersects $\sigma$ in the space $\tilde{X}$.
This means that $\iota'$ bijects a (vector space) basis of $R'_{n-i}$
with a basis of $\Hom_A({R}_{i}, A)$ and is therefore an isomorphism. A
similar argument shows that $\iota$ is an isomorphism as well.
\qed.

\begin{lemma}
There exists  an isomorphism of complexes
$$\tau_0:{N}\otimes_A{R}_*\to \Hom_A({R}'_{n-*},{N}).$$
\end{lemma}
\proof
We define $\tau_0$ as a composition,
$$\tau_0={m}\circ c\circ (id\otimes\iota): {N}\otimes_A{R}_*\to
\Hom_A({R}'_{n-*},{N}),$$
where $c$ is the canonical map 
$$c: N\otimes_A\Hom(R'_{n-*}, N)\to \Hom_A({R}'_{n-*},{N}\otimes A)$$ 
and $c: \Hom_A({R}'_{n-*},{N}\otimes A)\to\Hom_A({R}'_{n-*},{N})$ is the
map obtained from the action $A\otimes N\to N$.

It follows from a standard algebraic fact that when $R$ is a free
$A$-module of finite rank and $\iota$ is an isomorphism, the map $\tau_0$
is an isomorphism. (This is analogous to the fact that for
finite-dimensional vector spaces
$V$, $W$ there is an isomorphism $V^*\otimes W\to \Hom(V,W)$.)
\qed

This induces an isomorphism  $\tau:HH_*\to HH^{n-*}$ on homology.

Let 
$$\lambda=\tau\circ\kappa\circ\tau^{-1}: HH^*\to HH^{*-1}$$
 be 
the operator on $HH^*$ induced by the Connes operator $\kappa$ on $HH_*$
by $\tau$.

\begin{lemma}
The map $\lambda$ is a $BV$ operator compatible with the standard
Gerstenhaber algebra structure on $HH^*$.
\end{lemma}
\proof This follows from~\cite[Theorem 3.3.2]{gi2} (see
also~\cite[Sec.~6.5]{cbeg}). 
Ginzburg constructs a Poincar\'e duality isomorphism between $HH^i$ and
$HH_{n-i}$ of an algebra which has two
dual resolutions analogously to our construction of $\tau$. 
He shows using the formalism of non-commutative differential geometry
that this isomorphism sends the Connes operator to a BV operator on
$HH^*$. 
\qed

\

We have shown that the isomorphism $\rho$ takes the BV operator $\Delta$
on loop homology to the Connes operator $\kappa$ on $HH_*$ and that the
isomorphism $\tau$ takes $\kappa$ to a BV operator $\lambda$ on $HH^*$.
Let us define
$$\xi=\tau\circ\rho:\H_*\to HH^{n-*}.$$
It follows that $\xi$ is an isomorphism which takes $\Delta$ to
$\lambda$.

\begin{lemma}
The composition $$\xi=\tau\circ\rho:\H_*\to HH^{n-*}$$
is an associative algebra isomorphism.
\end{lemma}
\proof
The result follows from the following lemma.
\begin{lemma}
The map $\tau:HH_*\to HH^{n-*}$ takes the operation $\beta$ to cup
product on Hochschild cohomology.
\end{lemma}
\proof We first define cup product in terms of the complexes we have
used.

Let $_A Tot$ be the total complex of $R\otimes R'$ with the diagonal
$A$-action. 
\begin{lemma}
The complex $_ATot$ is an $A$-module resolution of $\mathbf{k}$.
\proof
This is a known result, see e.g.~\cite{w}.
\end{lemma}
We have a map $\cup_0: C^i\otimes C'^j\to \Hom_A(Tot_{i+j}, N)$
given by the multiplication on $N$ 
(which coincides with $A$ as a vector space). 
This means that the homology of the chain complex 
$\Hom_A(Tot_*,N)$ is the Hochschild cohomology, $HH^*(A)$.

\begin{lemma}
After passing to cohomology, the map $\cup_0$ becomes
the Hochschild cup product
$$
\cup:HH^{i}\otimes HH^j\to HH^{i+j}.
$$
\end{lemma}
\proof This is a standard fact from homological algebra,
see e.g.~\cite{loday}.
\qed

Using $\cup_0$ to compute the cup product, we show that this product is the
same as the product obtained from $\beta_0$.

\

We define a map of complexes $\beta'_0:\Hom_A(Tot_*, N)\to
A\otimes_{A\otimes A} W_{n-*}$ such that the following diagram commutes
after passing to homology
\begin{equation}
\label{simple-diagram}
\begin{diagram}
C_i \otimes C'_j & \rTo^{\beta_0} & A\otimes_{A\otimes A} W_{i+j-n}\\
\dTo^{\tau_0\otimes\tau_0} & & \uTo_\beta'_0 \\
C'^{n-i}\otimes C^{n-j} & \rTo^{\cup_0} & \Hom_A(Tot_{2n-i-j},N).
\end{diagram}
\end{equation}
Because the Poincar\'e duality map $\tau$ is obtained from an
augmentation of $\mu$, we see that $\beta'_0$ is a quasiisomorphism and
that on the level of homology, the corresponding map $\beta':HH^*\to
HH_{n-*}$ is the inverse of the isomorphism $\tau$.

This proves that the following diagram commutes.
\begin{equation}
\label{two-tier-diagram2}
\begin{diagram}
HH_i\otimes HH_j & \rTo^\beta & HH_{i+j-n}\\
\dTo^{\tau\otimes\tau} & & \dTo_\tau \\
HH^{n-i}\otimes HH^{n-j} & \rTo^\cup & HH^{2n-i-j} \quad \qed
\end{diagram}
\end{equation}

Combining this with our previous result, we obtain the commutative
diagram
\begin{equation}
\label{two-tier-diagram3}
\begin{diagram}
\H_i \otimes \H_j & \rTo^\bullet & \H_{i+j-n}\\
\dTo^{\rho\otimes\rho} & & \dTo_\rho \\
HH_i\otimes HH_j & \rTo^\beta & HH_{i+j-n}\\
\dTo^{\tau\otimes\tau} & & \dTo_\tau \\
HH^{n-i}\otimes HH^{n-j} & \rTo^\cup & HH^{2n-i-j}.
\end{diagram}
\end{equation}

\

Thus $\xi=\tau\circ\rho:\H_*\to HH_{n-*}$ is an
isomorphism of associative algebras.
We have also shown that it takes the BV operator
$\Delta$ to $\lambda$, the dual of the Connes operator $\kappa$.

A BV algebra is defined by its product and its BV operator, so we have
shown that the BV algebra structures on $\H_{n-*}$ and on $HH_*$ are
isomorphic. In particular, this implies that the Lie-like bracket on
string topology is mapped to the Gerstenhaber bracket on Hochschild
cohomology.
\samepage
This concludes the proof of Theorem \ref{thm:main}.
\quad Q.E.D.

% \pagebreak

\section{Concluding remarks}

\begin{enumerate}
\item 
 In this paper, we constructed an isomorphism of BV algebra
structures between the Hochschild cohomology and loop homology for
aspherical oriented closed manifolds. The algebra $\H_*(X)$ 
has additional algebraic structure (of a  2-dimensional
positive-boundary TQFT,  see~\cite{cg}) and it should be possible to
compute it algebraically for aspherical manifolds. 

\item  I hope that the methods of this work might be useful
for proving that the isomorphisms of algebras
$\H_{*-n}\cong HH^*(C(M),C(M))$ and
$\H_{*-n}\cong HH_*(C_*(\Omega(M)),C_*(\Omega(M)))$
for simply connected manifolds
also preserve BV structures.

\item 
I plan to extend the results of this paper for aspherical
orbifolds. String topology operations for orbifolds
have been recently introduced in~\cite{lux}.

\item 
It is known that the equivariant homology $H_*^{S^1}(LX)$
of $LX$ is related to the cyclic homology $HC_*(A)$ of the algebra
$A=C^*(X)$ of singular cochains of $X$ (see~\cite{jones}).  
It should be possible to show that for an aspherical manifold $X$ with
$\pi_1(X)=G$, structures on the cyclic cohomology $HC^*(\mathbf{k}[G])$ agree 
with the string topology operations on $H_*^{S^1}(LX)$.
\end{enumerate}

\

\noindent{\sc South Eugene High School, Eugene, OR}
\\[4pt]
{\tt mitkav@yahoo.com}
\end{document}